\documentclass[11pt]{article}
\usepackage{amsfonts}
\usepackage{mathrsfs}
\usepackage{amsmath}
\usepackage{amssymb}
\usepackage{graphicx}
\usepackage{epic}
\usepackage{amsthm,latexsym}
\usepackage{indentfirst}
\renewcommand{\paragraph}{\roman{paragraph}}
\usepackage[]{caption2}

\begin{document}

\title{\bf  Cacti with maximum Kirchhoff index}
\author{Wen-Rui Wang, Xiang-Feng Pan\thanks{Corresponding author. E-mail address: xfpan@ustc.edu.(X.-F.Pan), Ricciawang@163.com(W.-R.Wang).}, \\
{\small \it School of Mathematical Sciences, Anhui University,
Hefei 230601, P. R. China}}
\date{}
\maketitle

\noindent {\bf Abstract:}  The concept of resistance distance was first proposed by Klein and Randi\'c. The Kirchhoff index $Kf(G)$ of a graph $G$ is the sum of resistance distance between all pairs of vertices in $G$. A connected graph $G$ is called a cactus if each block of $G$ is either an edge or a cycle. Let $Cat(n;t)$ be the set of connected cacti possessing $n$ vertices and $t$ cycles, where $0\leq t \leq \lfloor\frac{n-1}{2}\rfloor$. In this paper, the maximum kirchhoff index of cacti are characterized, as well as the corresponding extremal graph.


\section{Introduction}
Topological indices and graph invariants based on the distances between the vertices of a graph are widely used in theoretical chemistry to establish relations between the structure and the properties of molecules. They provide correlations with physical, chemical and thermodynamic parameters of chemical compounds~\cite{Gutman2012Distance}.

All graphs considered in this paper will be finite, loopless, and contain no multiple edges. Let $G$ be a connected graph with vertex set $V(G)$, and edge set $E(G)$. The ordinary distance between $u$ and $v$, denoted by $d_G(u,v)$, is the number of edges in a shortest path joining $u$ and $v$ in $G$. For other undefined notations and terminology from graph theory, the readers are referred to~\cite{JA1976GraphTheory}. The famous Wiener index~\cite{Wiener1947} is equals to the sum of distances between all pairs of vertices, that is,
$$ W(G) =\sum_{i<j}d_{G}(v_{i},v_{j}).$$

On the basis of electrical network theory, Klein and Randi\'c~\cite{Klein1993Resistance} posed a new distance function named resistance distance in 1993. They viewed $G$ as an electrical network $N$ by replacing each edge of G with a unit resistor. The term resistance distance was used for the physical interpretation~\cite{Wang2010Cacti}: one imagines unit resistors on each edge of a connected graph $G$ with vertices $v_{1},v_{2},\cdots,v_{n}$ and takes the resistance distance between vertices $v_{i}$ and $v_{j}$ of $G$ to be the effective resistance between vertices $v_{i}$ and $v_{j}$, denoted by $R_G(v_{i},v_{j})$. This new kind of distance between vertices of a graph was diffusely studied in detail [1, 4, 6 - 12].

Analogous to the Wiener index, the Kirchhoff index~\cite{Klein1993Resistance} is defined as
$$ Kf(G) =\sum_{i<j}R_{G}(v_{i},v_{j}).$$
The Kirchhoff index of graphs as well as its application in chemistry attracts broad attention since it was put forward. Much work has been done to compute the Kirchhoff index of some classes of graphs, or give bounds for the Kirchhoff index of graphs and characterize extremal graphs. In unicyclic graphs extremal with respect to the Kirchhoff index were determined in~\cite{Yang2008Unicyclic,Zhang2009secondunicyclic}. Deng also studied the extremal Kirchhoff index of a class of unicyclic graphs~\cite{Guo2009unicyclic} and graphs with a given number of cut edges~\cite{Deng2010cutedges}. The extremal graphs with given matching number, connectivity, and minimal Kirchhoff index were characterized by Zhou in~\cite{Zhou2009matchingnumber}. Wang ~\cite{Wang2010Cacti} determined the first three minimal Kirchhoff indices among cacti. For further details and additional references, the readers may refer to [18 - 22].

Let G be a connected graph, and $deg_G(v)$ be the degree of a vertex $v$ in $G$. The longest path of $G$, considered in this paper, is a path with the largest resistance distance. A pendent vertex of G is a vertex of degree 1. Let $P_k = {r_1}{r_2} \cdots {r_k}(k\geq2)$ be a path of $G$ with distinct vertices $r_1, r_2,\cdots, r_k$ and assume that $deg_G(r_1) \geq 3, deg_G(r_2) = \cdots = deg_G(r_{k-1}) = 2,$ then $P_k$ is called a pendent path of length $k-1$ at $r_1$ in $G$ if $deg_G(r_k )=1$, and $P_k$ is called a internal path if $deg_G(r_k)\geq3$~\cite{Wang2015hyperWiener}. We suppose that $G_1$ and $G_2$ are two disjoint connected graphs with $u_1\in V(G_1)$ and $u_2 \in V(G_2)$. Let $(G,u) = (G_1,u_1) \oplus (G_2, u_2)$ denote the graph $G$ created by the coalescence of $u_1, u_2$, and denote the new common vertex by $u$. We refer to this procedure from $G_1$ and $G_2$ to $G$ as an identification operation~\cite{Wang2013Harary}.

A connected graph $G$ is called a $cactus$ if each block of $G$ is either an edge or a cycle. Denote by $Cat(n;t)$ the set of connected cacti possessing $n$ vertices and $t$ cycles, where $0\leq t \leq \lfloor\frac{n-1}{2}\rfloor$. We call $G$ a $star$ $cactus$ if all blocks share a common cut vertex. A $chain$  $cactus$ is a cactus graph $G$ if each block of it has at most two cut vertices and each cut vertex is shared by exactly two blocks. Any chain cactus with at least two blocks contains exactly two blocks that have only one cut-vertex. Such blocks are called $terminal$ $blocks$~\cite{Wang2015hyperWiener}. The cactus graph has many applications in real life problems, especially in radio communication system~\cite{Das2012algorithms}.

In this paper, the maximum kirchhoff index of cacti are characterized, as well as the corresponding extremal graph. The paper is organized as follows. In Section 2 we propose some useful graph operations changing the kirchhoff indices of graphs, we also cite some basic results used in the following sections. In Section 3 we find the extremal cacti with the maximum Kirchhoff and give explicit precise for the maximum kirchhoff index of cacti.

\section{Preliminaries}
\noindent {\bf Lemma 2.1}~\cite{Klein1993Resistance}  Let $x$ be a cut vertex of a connected graph $G$, and let $a$ and $b$ be vertices occurring in different components
which arise upon deletion of $x$. Then $R_G(a, b) = R_G(a, x) + R_G(x, b)$.

\noindent {\bf Lemma 2.2}~\cite{Wang2010Cacti}  Let $G_1$ and $G_2$ be connected graphs. If we identify any vertex, say $x_1$, of
$G_1$ with any other vertex, say $x_2$, of $G_2$ as a new common vertex $x$, and we obtain
a new graph $G$, then
$$Kf(G) = Kf(G_1) + Kf(G_2) + n_1 Kf_{x_2}(G_2) + n_2 Kf_{x_1}(G_1),$$
where $ Kf_{x_i}(G_i) =\sum_{{y_i}\in G_i} R_{G_i}(x_i,y_i),$ and $n_i =\mid V(G_i)\mid-1$ for $i=1,2$.

\noindent {\bf Lemma 2.3} Let $X$, $Y$ and $Z$  be connected graphs. Suppose that $a$ is a vertex of $Y$, $b$ is a vertex of $Z$, $x_1$ and $x_2$ are two endpoints of a longest path in $X$, $u$ and $v$ are two vertices of $X$ satisfying $0\leq R_X(x_1,u) \leq R_X(x_1,v) \leq R_X(x_1,x_2)$. Let $G_3$ be the graph constructed from $X$, $Y$, $Z$ by identifying $a$ with $u$ and $b$ with $v$, $G_4$ be formed from $X$, $Y$, $Z$ by identifying $a$ with $u$ and $b$ with $x_2$, $G_5$ be created from $X$, $Y$, $Z$ by identifying $a$ with $x_1$ and $b$ with $x_2$, as displayed in Fig. 1. Then

(1) when $Kf_v(X)\leq Kf_{x_2}(X)$, we have $Kf(G_3) < Kf(G_4)$;

(2) when $Kf_u(X)\leq Kf_{x_1}(X)$, we have $Kf(G_4) < Kf(G_5)$.

\noindent{\bf Proof.} (1) For convenience, we denote $x\in V(G)$ by $x\in G$. By the definition of kirchhoff index, one can get that

\begin{figure}
  \centering
  \includegraphics[width=5.5in]{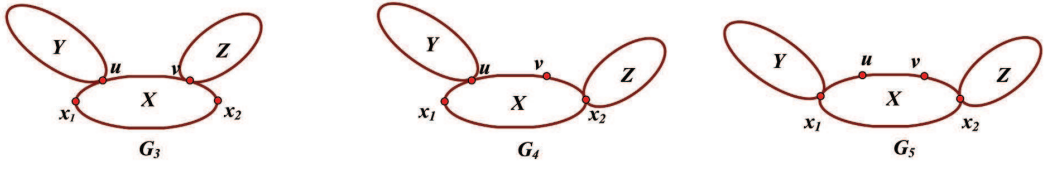}\\
  \caption{The graphs of Lemma 2.3}\label{fig1}
\end{figure}

\begin{eqnarray}
Kf(G_3) & = & \sum_{{x,y}\in Y-a} R_Y(x,y) + \sum_{{x,y}\in Z-b} R_Z(x,y) + \sum_{{x,y}\in X} R_X(x,y)
              + \underset{y\in Y-a}{\sum_{x\in X}}¡¡R_{G_3}(x,y) \nonumber\\
        & + & \underset{y\in Z-b}{\sum_{x\in X}}¡¡R_{G_3}(x,y) + \underset{y\in Z-b}{\sum_{x\in Y-a}}¡¡R_{G_3}(x,y) \nonumber
\end {eqnarray}
and
\begin{eqnarray}
Kf(G_4) & = & \sum_{{x,y}\in Y-a} R_Y(x,y) + \sum_{{x,y}\in Z-b} R_Z(x,y) + \sum_{{x,y}\in X} R_X(x,y)
              + \underset{y\in Y-a}{\sum_{x\in X}}¡¡R_{G_4}(x,y) \nonumber\\
        & + & \underset{y\in Z-b}{\sum_{x\in X}}¡¡R_{G_4}(x,y) + \underset{y\in Z-b}{\sum_{x\in Y-a}}¡¡R_{G_4}(x,y). \nonumber
\end {eqnarray}
Note that, for $x\in X, y\in Z-b,$
$$R_{G_4}(x,y)=R_Z(y,b)+R_X(x,x_2)\quad and \quad R_{G_3}(x,y)=R_Z(y,b)+R_X(x,v).$$
Then
$$R_{G_4}(x,y)-R_{G_3}(x,y)=R_X(x,x_2)- R_X(x,v).$$
For $x\in Y-a, y\in Z-b,$
$$R_{G_4}(x,y)-R_{G_3}(x,y)=R_X(u,x_2)- R_X(u,v)>0.$$
Hence, when $Kf_{x_2}(X)\geq Kf_v(X)$, it follows that

\begin{eqnarray}
Kf(G_4)-Kf(G_3) & = & \underset{y\in Z-b}{\sum_{x\in X}}(R_{G_4}(x,y) - R_{G_3}(x,y))
                       + \underset{y\in Z-b}{\sum_{x\in Y-a}}(R_{G_4}(x,y) - R_{G_3}(x,y)) \nonumber\\
                 & = & \underset{y\in Z-b}{\sum_{x\in X}}(R_X(x,{x_2}) - R_X(x,v))
                       + \underset{y\in Z-b}{\sum_{x\in Y-a}}(R_X(u,{x_2}) - R_X(u,v)) \nonumber\\
                 & > & (\mid Z \mid -1)(Kf_{x_2}(X)-Kf_v(X))\nonumber\\
                 & \geq & 0.\nonumber
\end {eqnarray}

(2) In a similar way to prove (1), we can prove (2). (For reference, please see  appendix.)

Then the proof of this lemma is complete.\hfill$\square$

\begin{figure}
  \centering
  \includegraphics[width=6.5in]{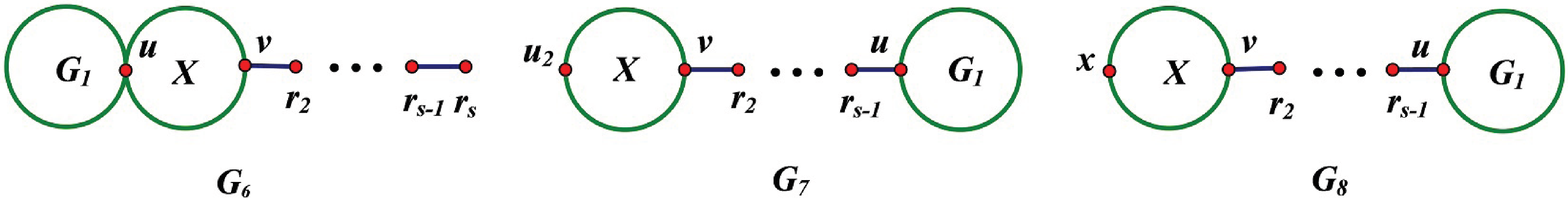}\\
   \includegraphics[width=3in]{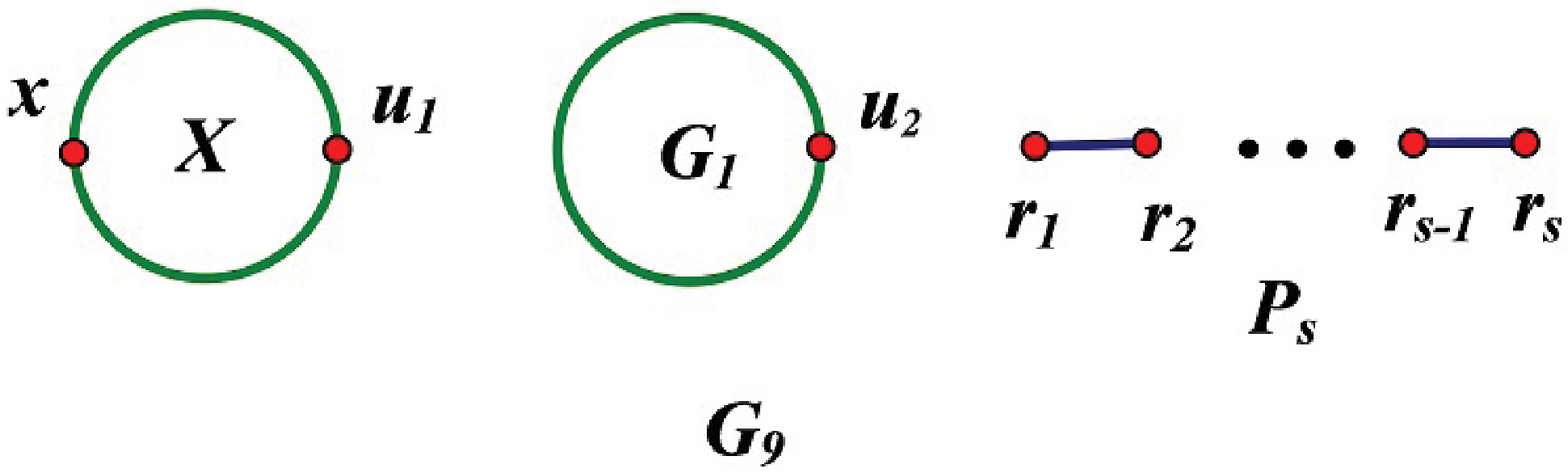}\\
  \caption{The graphs of Lemma 2.4}\label{fig2}
\end{figure}

The procedure of transformation from $G_3$ to $G_4$ is called an operation \uppercase\expandafter{\romannumeral 1}. Let $G_1$ and $X$ be two cacti not isomorphic to paths, $u_2$ and $x$ are two end vertices of a longest path in $X$. Then $G_6$ and $G_7$ are obtained by identification operation as follows:
\begin{eqnarray}
(M,v)  & = &(X,x)\oplus (P_s,r_1),\nonumber\\
(G_6,u)& = &(M,u_2)\oplus (G_1,u_1),\quad (G_7,u)=(M,r_s)\oplus (G_1,u_1).\nonumber
\end {eqnarray}

$G_8$ is constructed by attaching $G_1$ and $X$ to $r_1$, $r_s$ of $P_s$ , respectively. That is $(N,v) = (X,u_2)\oplus (P_s,r_1)$ and $(G_8,u)=(N,r_s)\oplus (G_1,u_1).$ Let $G_9$ be the group of $G_1$, $X$ and $P_s$. We call the process from $G_6$ to $G_7$ or $G_8$ an operation \uppercase\expandafter{\romannumeral 2}. $G_6$, $G_7$, $G_8$ and $G_9$ are drawn in Fig. 2.

\noindent {\bf Lemma 2.4}  Let $G_6$, $G_7$ and $G_8$ be three cacti defined above. For arbitrary connected graph $X$ containing at least one cycle, we have
$$Kf(G_6)< Kf(G_7) \quad or \quad Kf(G_6) < Kf(G_8).$$

\noindent{\bf Proof.} Firstly, we prove $Kf(G_6)< Kf(G_7)$.
By Lemma 2.2,
\begin{eqnarray}
Kf(G_6) & = & Kf(G_1)+Kf(M)+(\mid M \mid -1)Kf_{u_1}(G_1)+(\mid G_1 \mid -1)Kf_{u_2}(M),\nonumber\\
Kf(G_7) & = & Kf(G_1)+Kf(M)+(\mid M \mid -1)Kf_{u_1}(G_1)+(\mid G_1 \mid -1)Kf_{r_s}(M).\nonumber
\end {eqnarray}
It can be calculated that
\begin{eqnarray}
Kf_{u_2}(M) & = & Kf_{u_2}(X)+Kf_{u_2}(P_s)-R(u_2,x)=Kf_{u_2}(X)+(s-1)R(u_2,x)+Kf_{r_1}(P_s),\nonumber\\
Kf_{r_s}(M) & = & Kf_{r_s}(P_s)+Kf_{r_s}(X)-R(r_1,r_s)=Kf_{r_s}(P_s)+(s-1)(\mid {X} \mid -1)+Kf_x(X).\nonumber
\end {eqnarray}
Then we have
$$Kf_{r_s}(M)-Kf_{u_2}(M)=(s-1)(\mid X \mid -1 -R(u_2,x))+Kf_x(X)-Kf_{u_2}(X),$$
$$Kf(G_7)-Kf(G_6)=(\mid {G_1} \mid -1)((s-1)(\mid X \mid -1 -R(u_2,x))+Kf_x(X)-Kf_{u_2}(X)).$$
Note that
$$\mid {X} \mid > R(u_2,x) +1,$$
So we have $Kf(G_6) < Kf(G_7)$ when $Kf_x(X)\geq Kf_{u_2}(X).$
Similarly, we can get that $Kf(G_6)< Kf(G_8)$ when $Kf_{u_2}(X)\geq Kf_x(X).$ (For reference, please see  appendix.)
From the above arguments, the result follows.\hfill$\square$

\noindent {\bf Lemma 2.5}~\cite{Wang2010Cacti}  Let $u_1$ be a vertex of a connected graph $G_1$ and $w$, $x$ be two vertices
of a cycle $C_k (k > 3)$ with the largest $R(w, x)$. $G^{*}$ is a cactus depicted in Fig. 3, $(G^{**}, u) = (G_1, u_1) \oplus (C_k,w).$
Then  $Kf(G^{*})< Kf(G^{**}).$

As is shown in Fig. 3, $G_{10}$ is a $k$-vertex graph obtained by identifying one vertex from $C_3$ with one pendent vertex from $P_{k-2},$
while $G_{11}$ is obtained by identification operation :$(G_{11}, u) = (G_1, u_1)\oplus (G_{10}, r_{k-2}).$

\noindent {\bf Lemma 2.6} Let $G^{**}$ and $G_{11}$ be two cacti illustrated in Fig. 3. The procedure from
$G^{**}$ to $G_{11}$ is called an operation \uppercase\expandafter{\romannumeral 3}. Then $Kf(G^{**})< Kf(G_{11}).$

\noindent{\bf Proof.}  By Lemma 2.2,
\begin{eqnarray}
Kf(G^{**}) & = & Kf(G_1)+Kf(C_k)+(\mid C_k \mid -1)Kf_{u_1}(G_1)+(\mid G_1 \mid -1)Kf_w(C_k),\nonumber\\
Kf(G_{11}) & = & Kf(G_1)+Kf(G_{10})+(\mid G_{10} \mid -1)Kf_{u_1}(G_1)+(\mid G_1 \mid -1)Kf_{r_{k-2}}(G_{10}).\nonumber
\end {eqnarray}
Therefore,
\begin{eqnarray}
Kf(G_{11})-Kf(G^{**}) & = & Kf(G_{10})-Kf(C_k)+(\mid {G_1} \mid -1)(Kf_{r_{k-2}}(G_{10})-Kf_w(C_k))\nonumber\\
                      & + & (\mid G_{10} \mid - \mid C_k \mid )Kf_{u_1}(G_1).\nonumber
\end {eqnarray}
Also by Lemma 2.2, we have
$$Kf(G_{10})=Kf(C_3)+Kf(P_{k-2})+(\mid C_3 \mid -1)Kf_{r_1}(P_{k-2})+(\mid P_{k-2} \mid -1)Kf_{r_1}(C_3).$$
Recall that, $Kf(C_l)=\frac{l^3-l}{12}, \quad Kf_v(C_l)=\frac{l^2-1}{6}$ and $Kf(P_m)=\frac{m^3-m}{6}.$ Then one can have
\begin{eqnarray}
Kf(G_{10}) & = & \frac{1}{6}(k^3-11k+18),\nonumber\\
Kf_{r_{k-2}}(G_{10}) & = & \frac{1}{6}(3k^2-3k-10).\nonumber
\end {eqnarray}
Together with
\begin{eqnarray}
Kf(G_{10})-Kf(C_k) & = & \frac{1}{12}(k^3-21k+36),\nonumber\\
Kf_{r_{k-2}}(G_{10})-Kf_w(C_k) & = & \frac{1}{6}(2k^2-3k-9).\nonumber
\end {eqnarray}
Note that,
$$\mid G_{10} \mid = \mid C_k \mid \quad and \quad k > 3.$$
We arrive at
$$Kf(G_{11})-Kf(G^{**})=\frac{1}{12}(k-3)(k^2+3k-12)+\frac{1}{6}(\mid G_1 \mid -1)(2k+3)(k-3) > 0 . $$
Thus the conclusion holds.\hfill$\square$

\begin{figure}
  \centering
  \includegraphics[width=6.3in]{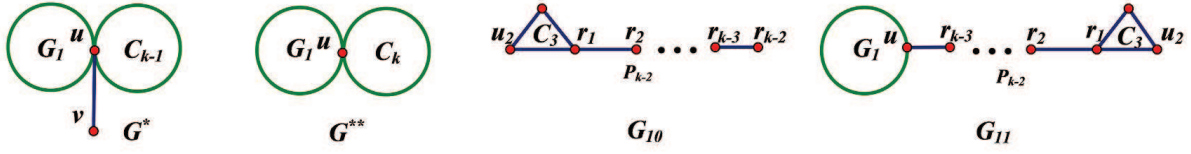}\\
  \caption{The graphs of Lemmas 2.5 and 2.6}\label{fig3}
\end{figure}

\begin{figure}
  \centering
  \includegraphics[width=4.2in]{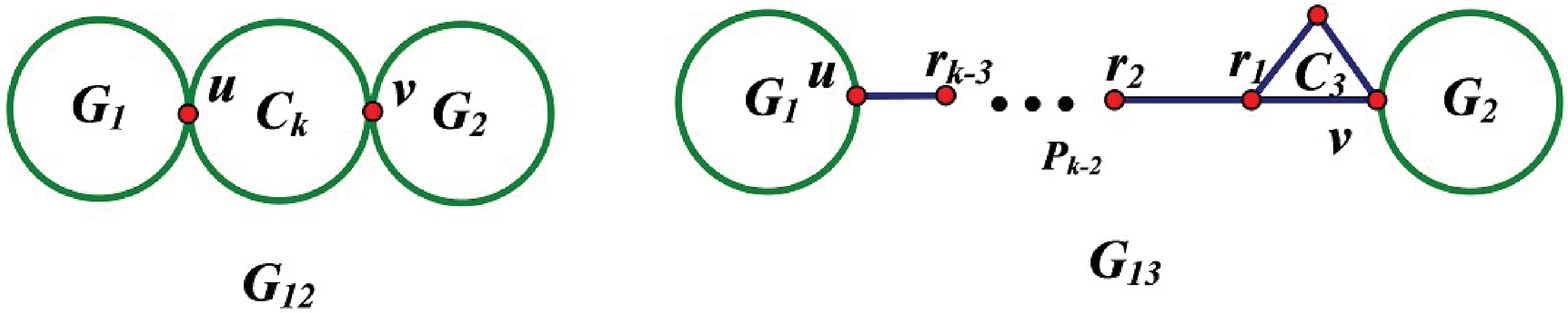}\\
  \caption{The graphs of Lemma 2.7}\label{fig4}
\end{figure}

\noindent {\bf Lemma 2.7}  Let $v_1$ be a vertex of connected graph $G_2$ and $u_2$ a vertex of $C_3$ in $G_{11}$ with $deg(u_2) = 2.$ $G_{12}$ and $G_{13}$ (see Fig. 4) are obtained by identification operation, namely
$$(G_{12}, v) = (G^{**}, x) \oplus (G_2, v_1),\quad (G_{13}, v) = (G_{11}, u_2) \oplus (G_2, v_1).$$
Operation \uppercase\expandafter{\romannumeral 4} is defined as the procedure from $G_{12}$ to $G_{13}.$ Then  $Kf(G_{12}) < Kf(G_{13}).$

\noindent{\bf Proof.}
By Lemma 2.2,
\begin{eqnarray}
Kf(G_{12}) & = & Kf(G^{**})+Kf(G_2)+(\mid G^{**} \mid -1)Kf_{v_1}(G_2)+(\mid G_2 \mid -1)Kf_x(G^{**}),\nonumber\\
Kf(G_{13}) & = & Kf(G_{11})+Kf(G_2)+(\mid G_{11} \mid -1)Kf_{v_1}(G_2)+(\mid G_2 \mid -1)Kf_{u_2}(G_{11}).\nonumber
\end {eqnarray}
Therefore, by Lemma 2.6
\begin{eqnarray}
Kf(G_{13})-Kf(G_{12}) & = & Kf(G_{11})-Kf(G^{**})+(\mid G_2 \mid -1)(Kf_{u_2}(G_{11})-Kf_x(G^{**})) \nonumber\\
                      & > & (\mid G_2 \mid -1)(Kf_{u_2}(G_{11})-Kf_x(G^{**})).\nonumber
\end {eqnarray}
It can be calculated that
\begin{eqnarray}
Kf_{u_2}(G_{11}) & = & Kf_{u_1}(G_1) + (\mid G_1 \mid -1)R(u,u_2) + Kf_{u_2}(G_{10}),\nonumber\\
Kf_x(G^{**}) & = & Kf_{u_1}(G_1) + (\mid G_1 \mid -1)R(u,x) + Kf_x(C_k).\nonumber
\end {eqnarray}

Then,
\begin{eqnarray}
Kf_{u_2}(G_{11})- Kf_x(G^{**}) & = & Kf_{u_2}(G_{10}) - Kf_x(C_k) + (\mid G_1 \mid -1)(R(u,u_2) - R(u,x)) \nonumber\\
                               & > & Kf_{u_2}(G_{10}) - Kf_x(C_k).\nonumber
\end {eqnarray}
Recall that, $Kf(C_l)=\frac{l^3-l}{12}, \quad Kf_v(C_l)=\frac{l^2-1}{6}$ and $Kf(P_m)=\frac{m^3-m}{6}.$ Then one can have
\begin{eqnarray}
Kf_{u_2}(G_{10}) & = & \frac{1}{6}(3k^2 - 11k + 14),\nonumber\\
Kf_x(C_k) &¡¡=¡¡&¡¡\frac{1}{6}(k^2-1).\nonumber
\end {eqnarray}
Therefore, since $k > 3$, it follows that
\begin{eqnarray}
Kf(G_{13})-Kf(G_{12}) & > & (\mid G_2 \mid -1)(Kf_{u_2}(G_{11})-Kf_x(G^{**})) \nonumber\\
                      & > & (Kf_{u_2}(G_{10}) - Kf_x(C_k)) \nonumber \\
                      & = &  \frac{1}{6}(2k-5)(k-3)\nonumber \\
                      & > &  0. \nonumber
\end {eqnarray}
This completes the proof.\hfill$\square$

\noindent {\bf Lemma 2.8}  Let $G_{13}$ be the graph defined in Lemma 2.7 and $G_{14}$ be the cactus of Fig. 5.
 Assume $a$ is the vertex of $C_3$ with $ deg_{G_{13}}(a) = 2$ in $G_{13}$. The process from $G_{14}$ to $G_{13}$ is referred to as an operation \uppercase\expandafter{\romannumeral 5}. Therefore, if $\mid G_1 \mid \leq \mid G_2 \mid ,$ then $Kf(G_{13}) \leq Kf(G_{14}) .$ Moreover, the equality
holds if and only if $\mid G_1 \mid = \mid G_2 \mid .$

\noindent{\bf Proof.}
Let $G_{10}^{'}$ and $G_{11}^{'}$ be connected graphs as indicated in Fig. 5.
It is easy to obtain that
\begin{eqnarray}
Kf(G_{13}) & = & Kf(G_{11})+Kf(G_2)+(\mid G_{11} \mid -1)Kf_{v_1}(G_2)+(\mid G_2 \mid -1)Kf_{u_2}(G_{11}),\nonumber \\
Kf(G_{14}) & = & Kf(G_{11}^{'})+Kf(G_1)+(\mid G_{11}^{'} \mid -1)Kf_{u_1}(G_1)+(\mid G_1 \mid -1)Kf_{r_{k-2}}(G_{11}^{'}).\nonumber
\end {eqnarray}
\begin{figure}
  \centering
  \includegraphics[width=2.7in]{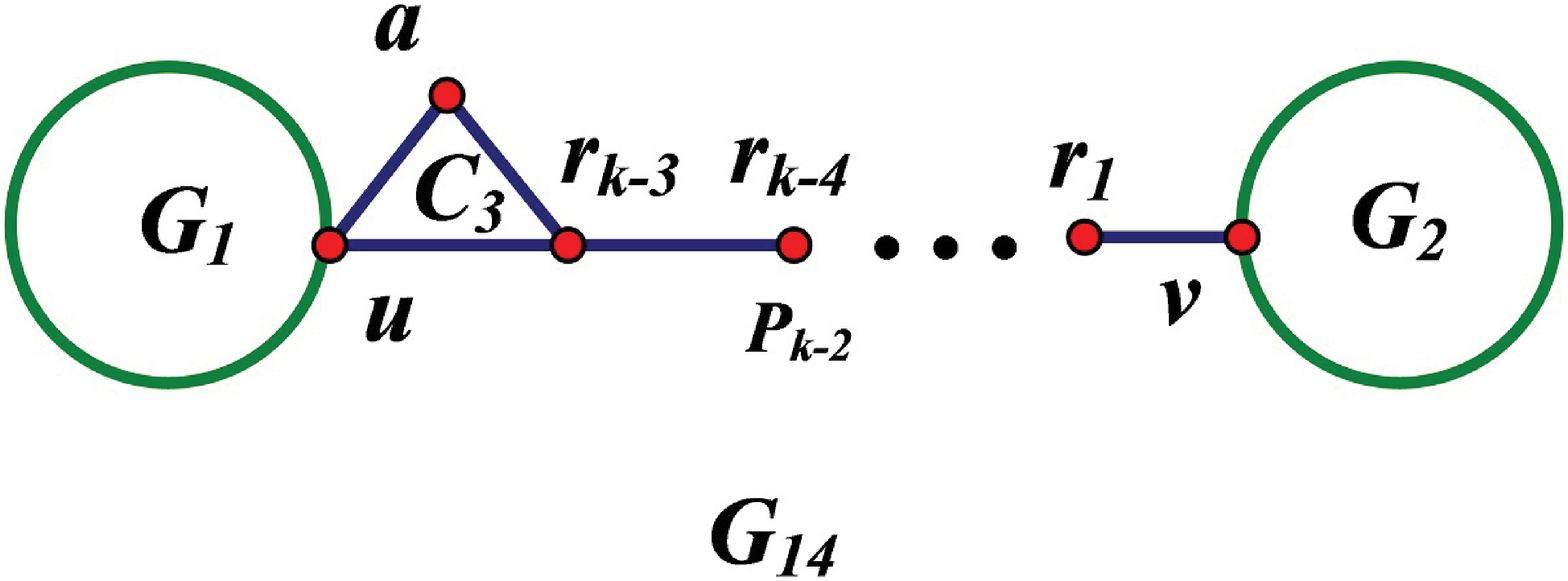}\\
   \includegraphics[width=4.5in]{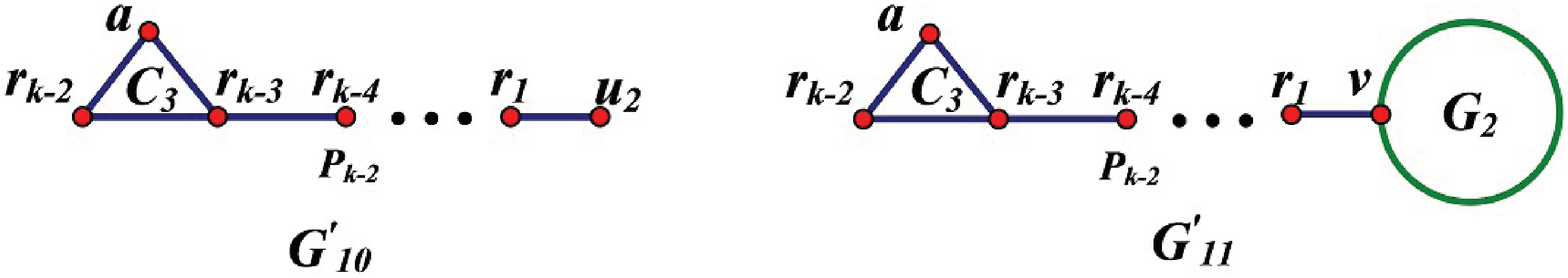}\\
  \caption{The graphs of Lemma 2.8}\label{fig5}
\end{figure}
Furthermore,
\begin{eqnarray}
Kf(G_{11}) & = & Kf(G_1)+Kf(G_{10})+(\mid G_{10} \mid -1)Kf_{u_1}(G_1)+(\mid G_1 \mid -1)Kf_{r_{k-2}}(G_{10}),\nonumber \\
Kf(G_{11}^{'}) & = & Kf(G_2)+Kf(G_{10}^{'})+(\mid G_{10}^{'} \mid -1)Kf_{v_1}(G_2)+(\mid G_2 \mid -1)Kf_{u_2}(G_{10}^{'}).\nonumber
\end {eqnarray}
Therefore,
\begin{eqnarray}
Kf(G_{14}) - Kf(G_{13}) & = & (\mid G_1 \mid -1)(Kf_{r_{k-2}}(G_{11}^{'}) - Kf_{r_{k-2}}(G_{10}))\nonumber \\
                        & + & (\mid G_2 \mid -1)(Kf_{u_2}(G_{10}^{'}) - Kf_{u_2}(G_{11}))\nonumber \\
                        & + & (\mid G_{11}^{'} \mid - \mid G_{10} \mid )Kf_{u_1}(G_1) + (\mid G_{10}^{'} \mid - \mid G_{11} \mid )Kf_{v_1}(G_2). \nonumber
\end {eqnarray}
Note that,
\begin{eqnarray}
Kf_{r_{k-2}}(G_{10}) & = & Kf_{u_2}(G_{10}^{'}) = \frac{1}{6}(3k^2-3k-10),\nonumber \\
Kf_{u_2}(G_{10}) & = & Kf_{r_{k-2}}(G_{10}^{'}) = \frac{1}{6}(3k^2 - 11k + 14)\nonumber
\end {eqnarray}
and
\begin{eqnarray}
Kf_{u_2}(G_{11}) & = & Kf_{u_1}(G_1) + (\mid G_1 \mid -1)R(u,u_2) + Kf_{u_2}(G_{10}),\nonumber \\
Kf_{r_{k-2}}(G_{11}^{'}) & = & Kf_{v_1}(G_2) + (\mid G_2 \mid -1)R(v_1,r_{k-2}) + Kf_{r_{k-2}}(G_{10}^{'}) .\nonumber
\end {eqnarray}
Thus,
\begin{eqnarray}
Kf(G_{14}) - Kf(G_{13}) & = & (\frac{4}{3}k - 4)(\mid G_2 \mid - \mid G_1 \mid )
                            -(\mid G_2 \mid + \mid G_{10} \mid - \mid G_{11}^{'} \mid -1)Kf_{u_1}(G_1)\nonumber \\
                        & + & (\mid G_1 \mid + \mid G_{10}^{'} \mid - \mid G_{11} \mid -1)Kf_{v_1}(G_2). \nonumber
\end {eqnarray}
Obviously,
\begin{eqnarray}
(\mid G_2 \mid + \mid G_{10} \mid - \mid G_{11}^{'} \mid -1) & = & (\mid G_2 \mid + \mid G_{10}^{'} \mid - \mid G_{11}^{'} \mid -1)=0,\nonumber \\
(\mid G_1 \mid + \mid G_{10}^{'} \mid - \mid G_{11} \mid -1) & = & (\mid G_1 \mid + \mid G_{10} \mid - \mid G_{11} \mid -1)=0.\nonumber
\end {eqnarray}
So we can get that
$$Kf(G_{14}) - Kf(G_{13}) = (\frac{4}{3}k - 4)(\mid G_2 \mid - \mid G_1 \mid ).$$
This implies that $Kf(G_{13}) \leq Kf(G_{14}) $ when $\mid G_1 \mid \leq \mid G_2 \mid $  and the equality holds if and only if $\mid G_1 \mid = \mid G_2 \mid .$
Thus the conclusion holds.\hfill$\square$

Lemma 2.3 to Lemma 2.8  indicate that the operations from \uppercase\expandafter{\romannumeral 1} to \uppercase\expandafter{\romannumeral 5} can lead to increasing the $Kf(G)$ of a cactus $G$.

\section{Main results}

In this section we characterize the extremal cacti with the largest Kirchhoff in $Cat(n;t)$. Assume that $Q_k$ is a chain cactus consisting of $k=\lfloor \frac{t}{2}\rfloor$ triangles. $u$ and $v$, satisfying $deg(u) = deg(v) = 2$, are two vertices from the terminal blocks of $Q_k$ and $Q_{t-k}$, respectively. Let $C_{n,t}$ (see Fig. 6) be another chain cactus obtained from $P_s$ , $Q_k$ and
$Q_{t-k}$ by identifying $u$ with $r_1$, $v$ with $r_s$. That is
$$(F,u)=(Q_k,u)\oplus (P_s,r_1),\quad (C_{n,t},v)=(F,r_s)\oplus (Q_{t-k},v).$$

\begin{figure}
  \centering
  \includegraphics[width=4.5in]{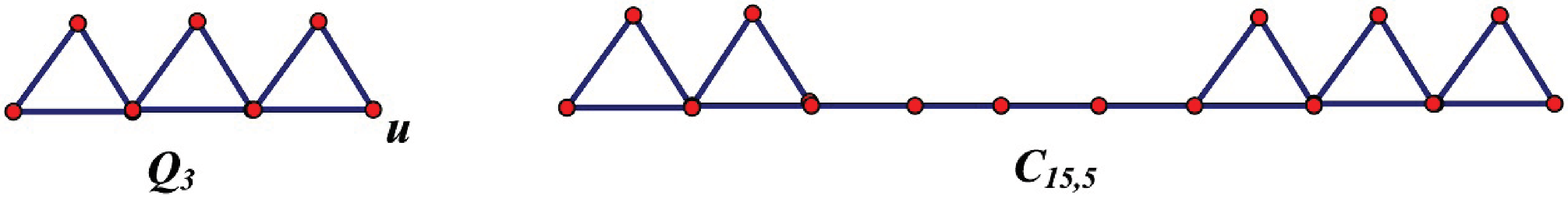}\\
    \includegraphics[width=3.7in]{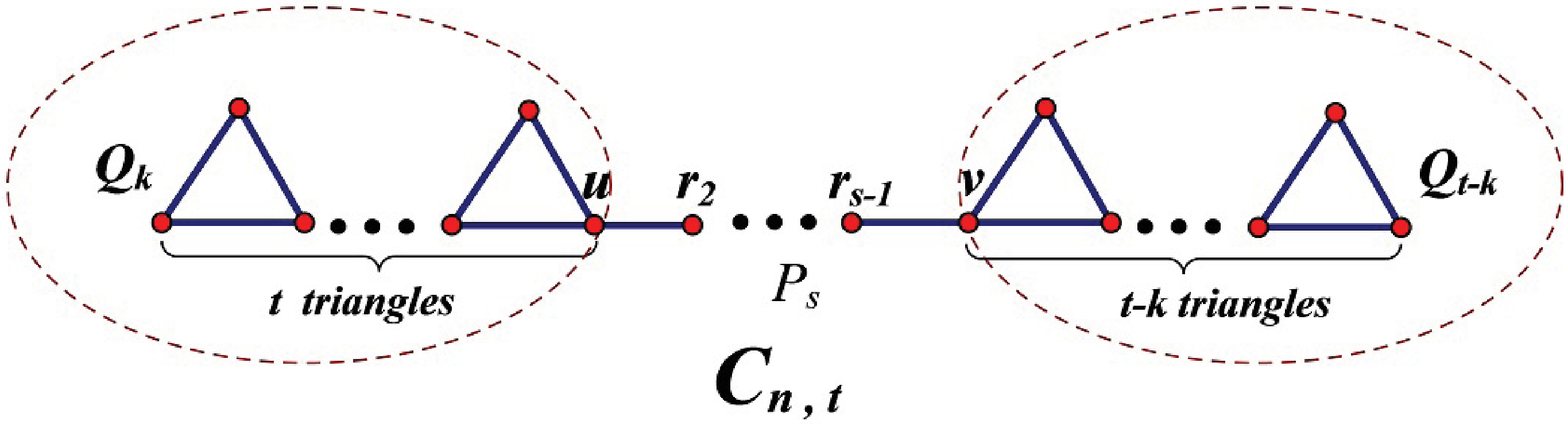}\\
  \caption{ The examples of $Q_k$ and $C_{n,t}$}\label{fig6}
\end{figure}

At first, we develop mathematical formulae describing $Kf$-values of $C_{n,t}$.
By the definitions of kirchhoff index, it is easy to check that
$$Kf_u(Q_k)=\frac{2}{3}(k^2+k).$$
Through accurate calculation, we arrive at
$$Kf_{r_s}(P_s)=\frac{1}{2}(s-1)s \quad and \quad Kf(P_s)=\frac{1}{6}(s^3-s).$$
While by Lemma2.2 and use of appropriate superposition computation, one can easy to obtain that
\begin{eqnarray}
Kf(Q_k) & = & Kf(Q_1)+Kf(Q_{k-1})+2Kf_u(Q_{k-1})+(\mid Q_{k-1} \mid -1)Kf_u(Q_1) \nonumber \\
        & = & \frac{2}{9}(2k^3+6k^2+k). \nonumber
\end {eqnarray}

Similarly, for $F$, we can get that
\begin{eqnarray}
Kf_{r_s}(F) & = & Kf_{r_s}(P_s)+Kf_u(Q_k)+(\mid Q_k \mid -1)(s-1) \nonumber \\
            & = & \frac{1}{2}(s^2-s) + \frac{2}{3}(k^2+k) + 2k(s-1) \nonumber
\end {eqnarray}
and
\begin{eqnarray}
Kf(F) & = & Kf(Q_k)+Kf(P_s)+(\mid Q_k \mid -1)Kf_{r_s}(P_s)+(\mid P_s \mid -1)Kf_u(Q_k) \nonumber \\
      & = & \frac{2}{9}(2k^3+3k^2-2k) + \frac{1}{6}(s^3-s) + \frac{1}{3}ks(2k+3s-1). \nonumber
\end {eqnarray}
So, for $C_{n,t}$, we have
$$Kf(C_{n,t})= Kf(F) + Kf(Q_{t-k})+(\mid F \mid -1)Kf_v(Q_{t-k})+(\mid Q_{t-k} \mid -1)Kf_{r_s}(F).$$

Note that, $n=2t+s$ and $ k = \lfloor \frac{t}{2} \rfloor = \left\{\begin{array}{ll}\frac{t}{2}, &\mbox{ if }\text{t is even;}\\
                                                 \frac{t-1}{2}, &\mbox{ if }\text{t is odd.}
                                 \end{array}\right.$

Thus, by the definition of kirchhoff index and suitable precise computation, it can be calculated that

$ Kf(C_{n,t}) = \left\{\begin{array}{ll}\frac{1}{18}(3n^3-3n-12nt^2-6nt+8t^3+12t^2-2t), &\mbox{ if }\text{t is even;}\\
                                                 \frac{1}{18}(3n^3-15n-12nt^2-6nt+8t^3+12t^2+22t+12), &\mbox{ if }\text{t is odd.}
                                 \end{array}\right.$

Next, we demonstrate that $C_{n,t}$ is the extremal cacti with the maximum $Kf$. For $Cat(n;t)$, $Cat(n;0)$ is the set of all trees and $Cat(n;1)$ is the set of unicyclic graphs. Their extremal graphs are $Cat(n;0) = P_n$ and $Cat(n;1)$, respectively. So we suppose that $n\geq 5$ and $t \geq 2$.

\noindent {\bf Theorem 3.1}
$G \in Cat(n; t ) - \{C_{n,t} \}$ for $n \geq 5$ and $t \geq 2$, then
 $$Kf(G)<Kf(C_{n,t}).$$
\noindent{\bf Proof.}  The proof contains the following four steps.

Firstly, suppose that $G$ is a connected graph, but not isomorphic to chain cactus, then a chain cactus $H_1$ is created by repeated applications of the operation  \uppercase\expandafter{\romannumeral 1} described in Lemma 2.3. By Lemma 2.3, it follows that
 $$Kf(G) < Kf(H_1).$$
Obviously, all the cut vertices of $H_1$ are in a longest path.

In the second step, our goal is convert pendent paths into internal paths.
At present, $H_1$ consists of at most two pendent paths, which are only possibly  fused to the first and last cycles, respectively. Assume that $H_1$ has one pendent path $P_t$ or two. Without loss of generality, one can suppose that $P_t$ is attached to the first cycle in $H_1$, then $P_t$ can be transformed into internal path by operation \uppercase\expandafter{\romannumeral 2} defined in Lemma 2.4. We denote the new graph by $H_1^{'}$. From Lemma 2.4, one can easily conclude that
$$Kf(H_1) < Kf(H_1^{'}).$$
In a similar way, we can have the same results if there is another pendent path attached to the last cycle. Thus, there is no pendent path in the newly constructed graph $H_2$.

Now, it comes to the third step. In this step, we will use two operations mentioned above.
Let $C_l$  $(l > 3)$ be any cycle of $H_2$. If $C_l$ contains only one cut vertex, then use operation \uppercase\expandafter{\romannumeral 3} established in Lemma 2.6 ; If $C_l$ contains two cut vertices, then apply operation \uppercase\expandafter{\romannumeral 4} described in Lemma 2.7 to it.
After finite times transformations in this step, we rapidly construct a cactus $H_3$, in which every cycle is triangle. According to Lemmas 2.5, 2.6 and 2.7, we arrive at
$$Kf(H_2)<Kf(H_3).$$

Finally, let $H_3 \ncong C_{n,t}$. Suppose that $H_3$ contains only one internal path which is not between $Q_k$ and $Q_{t-k}$. We use operation \uppercase\expandafter{\romannumeral 5} defined in Lemma 2.8 for appropriate times, then we have the chain cactus $C_{n,t}$. If there are more than one internal path in $H_3$, fuse every internal path into a unique internal path connecting $Q_k$ and $Q_{t-k}$ by operation \uppercase\expandafter{\romannumeral 5}, then $C_{n,t}$ is also created. we have $$Kf(H_3)<Kf(C_{n,t})$$ based on Lemma 2.8.

Therefore, by the four steps above, we constructed chain cactus $C_{n,t}$. Since $G$ is not isomorphic to chain cactus, it is not difficult to demonstrate that $G \ncong C_{n,t}$, thus
$$Kf(G)< Kf(C_{n,t}).$$
\hfill$\square$

\begin{figure}
  \centering
  \includegraphics[width=2in]{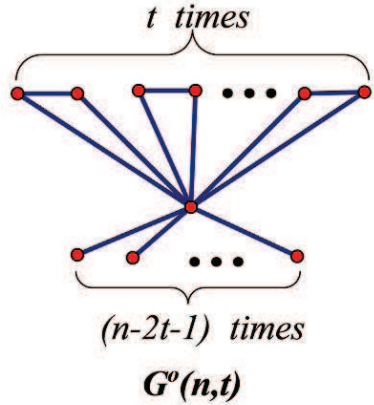}\\
  \caption{Cactus with the minimum Kirchhoff index }\label{fig7}
\end{figure}

In the end of this section, we give upper and lower bounds of $Kf$-value of a cactus. In~\cite{Wang2010Cacti}, the minimum Kirchhoff index and the corresponding extremal graph, denote by $G^0(n,t)$ (see Fig. 7), are presented by Wang and Hua.  Combining the results given by them, we arrive at the following proposition.

\noindent {\bf Proposition 3.2} Let $G$ be any $n$-vertex cactus different from $C_{n,t}$ and $G^0(n, t )$, then
$$Kf(G^0(n, t )) < Kf(G) < Kf(C_{n,t}).$$
\\
\\
\\

\begin{eqnarray}
Kf(G_5) & = & \sum_{{x,y}\in Y-a} R_Y(x,y) + \sum_{{x,y}\in Z-b} R_Z(x,y) + \sum_{{x,y}\in X} R_X(x,y)
              + \underset{y\in Y-a}{\sum_{x\in X}}¡¡R_{G_5}(x,y) \nonumber\\
        & + & \underset{y\in Z-b}{\sum_{x\in X}}¡¡R_{G_5}(x,y) + \underset{y\in Z-b}{\sum_{x\in Y-a}}¡¡R_{G_5}(x,y), \nonumber
\end {eqnarray}
and
\begin{eqnarray}
Kf(G_4) & = & \sum_{{x,y}\in Y-a} R_Y(x,y) + \sum_{{x,y}\in Z-b} R_Z(x,y) + \sum_{{x,y}\in X} R_X(x,y)
              + \underset{y\in Y-a}{\sum_{x\in X}}¡¡R_{G_4}(x,y) \nonumber\\
        & + & \underset{y\in Z-b}{\sum_{x\in X}}¡¡R_{G_4}(x,y) + \underset{y\in Z-b}{\sum_{x\in Y-a}}¡¡R_{G_4}(x,y). \nonumber
\end {eqnarray}
Note that, for $x\in X, y\in Y-a,$
$$R_{G_5}(x,y)=R_Y(y,a)+R_X(x,x_1)\quad and \quad R_{G_4}(x,y)=R_Y(y,a)+R_X(x,u).$$
Then
$$R_{G_5}(x,y)-R_{G_4}(x,y)=R_X(x,x_1)- R_X(x,u).$$
For $x\in Y-a, y\in Z-b,$
$$R_{G_5}(x,y)-R_{G_4}(x,y)=R_X(x_1,x_2)- R_X(u,x_2)>0.$$
Hence, when $Kf_{x_1}(X)\geq Kf_u(X)$, it follows that

\begin{eqnarray}
Kf(G_5)-Kf(G_4) & = & \underset{y\in Y-a}{\sum_{x\in X}}(R_{G_5}(x,y) - R_{G_4}(x,y))
                       + \underset{y\in Z-b}{\sum_{x\in Y-a}}(R_{G_5}(x,y) - R_{G_4}(x,y)) \nonumber\\
                 & = & \underset{y\in Y-a}{\sum_{x\in X}}(R_X(x,{x_1}) - R_X(x,u))
                       + \underset{y\in Z-b}{\sum_{x\in Y-a}}(R_X({x_1},{x_2}) - R_X(u,x_2)) \nonumber\\
                 & > & (\mid Y \mid -1)(Kf_{x_1}(X)-Kf_u(X))\nonumber\\
                 & \geq & 0.\nonumber
\end {eqnarray}
\hfill$\square$

\noindent{\bf Proof of Lemma 2.4(2).}
  Here, we prove $Kf(G_6)< Kf(G_8)$.
By Lemma 2.2,
\begin{eqnarray}
Kf(G_6) & = & Kf(G_1)+Kf(M)+(\mid M \mid -1)Kf_{u_1}(G_1)+(\mid G_1 \mid -1)Kf_{u_2}(M),\nonumber\\
Kf(G_8) & = & Kf(G_1)+Kf(N)+(\mid N \mid -1)Kf_{u_1}(G_1)+(\mid G_1 \mid -1)Kf_{r_s}(N).\nonumber
\end {eqnarray}
It can be calculated that
\begin{eqnarray}
Kf(M) & = & Kf(X)+Kf(P_s)+(s -1)Kf_x(X)+(\mid X \mid -1)Kf_{r_1}(P_s),\nonumber\\
Kf(N) & = & Kf(X)+Kf(P_s)+(s -1)Kf_{u_2}(X)+(\mid X \mid -1)Kf_{r_1}(P_s)\nonumber
\end {eqnarray}
and
\begin{eqnarray}
Kf_{u_2}(M) & = & Kf_{u_2}(X)+Kf_{u_2}(P_s)-R(u_2,x)=Kf_{u_2}(X)+(s-1)R(u_2,x)+Kf_{r_1}(P_s);\nonumber\\
Kf_{r_s}(N) & = & Kf_{r_s}(P_s)+Kf_{r_s}(X)-R(r_1,r_s)=Kf_{u_2}(X)+(s-1)(\mid {X} \mid -1)+Kf_{r_s}(P_s).\nonumber
\end {eqnarray}
Note that
$$\mid {X} \mid > R(u_2,x) +1,$$
Then we have
$$Kf_{r_s}(N)-Kf_{u_2}(M) = (s-1)(\mid X \mid -1 -R(u_2,x)) >0 .$$
Hence
\begin{eqnarray}
Kf(G_8)-Kf(G_6) & = & Kf(N)-Kf(M) + (\mid {G_1} \mid -1)(Kf_{r_s}(N)-Kf_{u_2}(M))\nonumber\\
                & = & (s-1)(Kf_{u_2}(X)-Kf_x(X)) + (\mid {G_1} \mid -1)(Kf_{r_s}(N)-Kf_{u_2}(M))\nonumber\\
                & > & (s-1)(Kf_{u_2}(X)-Kf_x(X)).\nonumber
\end {eqnarray}

So we have $Kf(G_6)< Kf(G_8)$ when $Kf_{u_2}(X)\geq Kf_x(X).$ \hfill$\square$

\end{document}